\documentclass[11pt]{article}
\usepackage[round]{natbib}
\usepackage{mathrsfs}
\usepackage{latexsym}
\usepackage{amsmath,amsthm}
\usepackage{amsfonts}
\usepackage{url}

\usepackage{amssymb}

\newcommand{\g}{\gamma}

\usepackage{bbm}
\newcommand{\chr}{\boldsymbol{\mathbbm{1}}} % characteristic function
\newcommand{\pred}[1]{\chr_{\left\{ #1 \right\}}}
\newcommand{\prs}{\vec{P}}
\newcommand{\pr}[1]{\prs\!\tlprn{#1}}
\newcommand{\mexp}{\vec{E}}

\newcommand{\E}{\mexp}
\renewcommand{\P}{\prs}
\newcommand{\Pn}{\P\!_n}

\newcommand{\Dn}{\Delta_n}

\newcommand{\Given}{\, \labs \,}
\newcommand{\etab}{\bar\eta}
\newcommand{\TV}[1]{\nrm{#1}_{\textrm{{\tiny \textup{TV}}}}}
\newcommand{\ham}{\textrm{{\tiny \textup{Ham}}}}
\newcommand{\Lip}[1]{\nrm{#1}_{\textrm{{\tiny \textup{Lip}}}}}
\DeclareMathOperator{\supp}{supp}

\newcommand{\convd}{\Rightarrow}

\newcommand{\diam}{\operatorname{diam}}
\newcommand{\essup}{\mathop{\operatorname{ess\,sup}}}

\newcommand{\Lone}{L_1(\R^n,\mu^n)}

\newcommand{\ohs}[2]{
\!
\begin{array}{ll}
{#1} \\
{#2}
\end{array}
\!
}
\newcommand{\oh}[2]{
\!
\begin{array}{ll}
\observed{#1} \\
\hidden{#2}
\end{array}
\!
}
\newcommand{\hp}[1]{{\hidden{#1}}'}
\newcommand{\op}[1]{{\observed{#1}}'}

\newcommand{\subplus}{_{
+
}}

\newcommand{\pl}[1]{\paren{#1}\subplus}

\usepackage{geometry}
\makeatletter
\theoremstyle{plain}
\newtheorem{thm}{Theorem}[section]
\theoremstyle{plain}
\newtheorem*{thm*}{Theorem}
 \newtheorem{cor}[thm]{Corollary} %%Delete [thm] to re-start numbering
 \theoremstyle{plain}    
 \theoremstyle{remark}
 \newtheorem{rem}[thm]{Remark}
 \theoremstyle{plain}
 \newtheorem{lem}[thm]{Lemma} %%Delete [thm] to re-start numbering
\makeatother
\newcommand{\bethn}{\begin{thm}}
\newcommand{\enthn}{\end{thm}}
\renewcommand{\beth}{\begin{thm*}}
\newcommand{\enth}{\end{thm*}}
\newcommand{\bepf}{\begin{proof}}
\newcommand{\enpf}{\end{proof}}
\newcommand{\belen}{\begin{lem}}
\newcommand{\enlen}{\end{lem}}
\newcommand{\ben}{\begin{enumerate}}
\newcommand{\een}{\end{enumerate}}
\newcommand{\bit}{\begin{itemize}}
\newcommand{\eit}{\end{itemize}}
\newcommand{\becon}{\begin{cor}}
\newcommand{\encon}{\end{cor}}

\newcommand{\hsymb}{\mathrm{h}}
\newcommand{\osymb}{\mathrm{o}}
\newtheorem{assumption}[thm]{Assumption}
\newcommand{\be}{\begin{equation}}
\newcommand{\ee}{\end{equation}}
\newcommand{\hiddenX}{\Omega_{\hsymb}}
\newcommand{\observedX}{\Omega_{\osymb}}
\newcommand{\observedx}{{x^{\osymb}}}
\newcommand{\hiddenx}{{x^{\hsymb}}}
\newcommand{\observed}[1]
{#1^{\osymb}}
\newcommand{\hidden}[1]
{#1^{\hsymb}}
\newcommand{\h}[1]{\hidden{#1}}
\renewcommand{\o}[1]{\observed{#1}}

\newtheorem{definition}[thm]{Definition}

\renewcommand{\vec}[1]{\bs{\mathrm{#1}}}

\newcommand{\phinorm}[1]{\nrm{#1}_{\Phi}}
\newcommand{\psinorm}[1]{\nrm{#1}_{\Psi}}
\newcommand{\basicspace}{\Omega}
\newcommand{\X}{\basicspace}

\newcommand{\seq}[3]{#1_{#2}\ldots#1_{#3}}
\newcommand{\sseq}[3]{#1_{#2:#3}}  % short seq

\newcommand{\sd}{,}
\newcommand{\scat}[1]{(#1)}

\newcommand{\labs}{\left| \vphantom{\sum_a^b} \right.}

\newcommand{\nrm}[1]{\left\Vert #1 \right\Vert}

\newcommand{\tsnrm}[1]{\Vert #1 \Vert}
\newcommand{\iprod}[1]{\left\langle #1 \right\rangle}

\newcommand{\phix}[1]{\Phi_{#1}}
\newcommand{\psix}[1]{\Psi_{#1}}
\newcommand{\iprodx}[2]{\left\langle #2 \right\rangle_{#1}}

\newcommand{\f}{\varphi}
\renewcommand{\k}{\kappa}

\newcommand{\psin}[1]{\Psi_{#1}}

\newcommand{\calA}{\mathcal{A}}

\newcommand{\calF}{\mathcal{F}}

\newcommand{\calL}{\mathcal{L}}

\newcommand{\calP}{\mathcal{P}}

\newcommand{\calW}{\mathcal{W}}
\newcommand{\calX}{\mathcal{X}}
\newcommand{\calY}{\mathcal{Y}}

\newcommand{\tha}{\theta}
\newcommand{\R}{\mathbb{R}}
\newcommand{\N}{\mathbb{N}}

\newcommand{\QQ}{\mathbb{Q}}
\newcommand{\beq}{\begin{eqnarray*}}
\newcommand{\eeq}{\end{eqnarray*}}
\newcommand{\beqn}{\begin{eqnarray}}
\newcommand{\eeqn}{\end{eqnarray}}
\newcommand{\paren}[1]{\left( #1 \right)}

\newcommand{\tlprn}[1]{\left\{ #1 \right\}}
\newcommand{\set}[1]{\tlprn{#1}}
\newcommand{\abs}[1]{\left| #1 \right|}
\newcommand{\ceil}[1]{\ensuremath{\left\lceil#1\right\rceil}}
\newcommand{\gn}{\, | \,}
\newcommand{\ts}{\textstyle}
\newcommand{\bs}{\boldsymbol}
\renewcommand{\th}{\ensuremath{^{\mathrm{th}}}~}

\newcommand{\hide}[1]{}
\newcommand{\oo}[1]{\frac{1}{#1}}

\def\eps{\varepsilon}
\newcommand{\defeq}{=}
\newcommand{\abscont}{\ll}
\newcommand{\sigalg}{\calF}
\newcommand{\tp}{\otimes}

\title{A Strong Law of Large Numbers for
Strongly Mixing Processes
}
\author{L. Kontorovich and A. Brockwell}

\begin{document}
\maketitle
\begin{abstract}
We prove a strong law of large numbers for a class of strongly mixing
processes.
%The 
Our
result rests on
recent advances in understanding of concentration of measure.
It is simple to apply and gives finite-sample (as
opposed to asymptotic) bounds, with readily computable
rate constants.  In particular, this makes it
suitable for analysis of inhomogeneous Markov processes.
We demonstrate how it can be applied to establish an
almost-sure convergence result
for a class of models that includes
as a special case a class of adaptive Markov chain Monte Carlo
algorithms.
\end{abstract}

\section{Introduction}

The strong laws of large numbers (SLLNs) play a fundamental
role in statistics.  They assert the
convergence of empirical averages to true expectations, and, under
appropriate assumptions, ensure that inferences
about persistent world phenomena
become increasingly more valid
%behave in a desirable way
as data accumulates.
The various laws of large numbers date back at least to
the publication in 1713 of Jakob Bernoulli's
\emph{Ars Conjectandi} \citep{bernoulli1713},
which stated an early form
of the weak law of large numbers.
Subsequent development of the concept was carried out
by many others.  Among the first of these was
Sim\'{e}on-Denis Poisson, who generalized
Bernoulli's result and gave it its modern name,
``la loi de grands nombres".  Other notable mathematicians
who made further contributions include Chebyshev,
Markov, Borel, Cantelli, and Kolmogorov.
Over the years, the various contributions gave rise to
two common forms, the weak and strong laws of large numbers,
establishing conditions, respectively, for weak and strong convergence
of empirical averages.

The SLLN in its earlier forms
applies to sequences of independent random variables.
However, for dependent sequences the theory is not as well
developed.  Indeed, to quote \cite{ninness00},
``for non-iid sequences
the required SLLN results do not seem to be readily
available in the literature''
and thus researchers continue to work to
develop generalized strong laws of large numbers.
Since around 1960, a number 
of 
different generalizations have been obtained.
%people have
%generalized it in various ways.
Our goal is not to give a comprehensive survey of these results;
the interested reader may find
%but
a useful list of relevant papers
%can be found
at
{\tt http://www.stats.org.uk/law-of-large-numbers/}.
%It may be presumptuous to attempt to condense a large body of work to
%a pithy description
%(especially since we have not read every single one of the papers on
%the topic).
%Not surprisingly,
%much of this work has coincided with the formalization
%of the various kinds of mixing properties that processes may exhibit.
The basic current state of affairs may be
summarized as follows. From Birkhoff's ergodic theorem we get a law of
large numbers for ergodic processes; this has been strengthened by
\cite{breiman60} to cover the case where the stationary
distribution is singular with respect to the Lebesgue
measure. Assumptions of ergodicity are typically too weak to provide a
convergence rate -- this requires a stronger mixing condition. A
classical (and perhaps first of its kind) example of the latter is the
paper by
%Blum, Hanson and Koopmans
\cite{blum63}, which proves a
strong law of large numbers under a mixing condition known in a modern
form as $\psi$-mixing (see \cite{bradley05}).
This mixing condition guarantees exponentially rapid convergence,
but the proof does not directly yield rate constants.

We approach the problem of developing such a SLLN
with easily computable rate constants.  This is useful
for example, in the context of simulation-based algorithms,
 if one is interested in determining how many
iterations are required to achieve a specified accuracy.
To this end, we assume a stronger (though still quite realistic) mixing
condition from which we
%Of the recent work, one of the most general forms of
%the SLLN for dependent processes is stated in
%\cite{blum63}.
%{\bf probably also Breiman paper}
%It is also worth noting that
%many SLLNs give asymptotic rates but not finite sample bounds. For
%many practical applications, one actually wants to know for how many
%steps to run an algorithm to achieve a specified accuracy at a given
%confidence level.
%In this paper, we generalize the results of
%\cite{blum63}, and provide explicit finite sample
%bounds, making use of concentration of measure results.
%Using this approach we shall be able to
obtain strong laws of large numbers, along with finite-sample bounds,
for sequences of random variables with arbitrary dependence.
In a sense, the application of concentration of measure
theory to establish a SLLN is straightforward, although
an important technical contribution of this paper is to extend the
deviation bounds proved in
\cite{kontram06,kont06-hmm,kont06-metric} for discrete measure
spaces to the continuous case,
as well as to bound certain mixing coefficients for adaptive Markov
chains
and to prove a strong law of large numbers for these (Theorem \ref{thm:adaptslln}).
We will
state our results in a fairly general setting.
We emphasize that the
basic
concentration result presented and employed
here (Theorem \ref{thm:mgalebd})
%are
is
by no means the only one available.
The recent result of
\cite{chazottes07} yields a very similar inequality; earlier results
along these lines appeared in
\cite{marton98}.\hide{
Marton's inequality
\cite{marton95} could be converted into a requisite law of large
numbers for the more restrictive class of contracting Markov chains; a
more recent result of
%Chazotte et al.
\cite{chazottes07} yields an
inequality of a very similar type to what we shall be applying. }
The
interested reader may refer to authoritative and comprehensive surveys
of concentration techniques, such as
\cite{ledoux01},
\cite{schechtman99}, and
\cite{lugosi03}
or to \cite{chatterjee05} and
\cite{kont07-thesis} for more recent developments.

We mention in passing that the strong laws we consider here are not
{\em uniform laws}, the latter asserting almost sure convergence
uniformly over some permissible (i.e., Glivenko-Cantelli\footnote{
Following \cite{pollard84},
the term {\em permissible} is used to avoid measure-theoretic
pathologies associated with taking suprema over uncountable
collections of sets.}) class of sets. There has been some work on such
uniform laws for non-i.i.d. processes. In particular, \cite{nobel93} show
that if $\calF$ is a Glivenko-Cantelli class for an i.i.d. process,
then the corresponding uniform strong law also holds for
$\beta$-mixing
(absolutely regular)
%(weakly Bernoulli)
process. In the other direction,
\cite{nobel95} gives a counterexample where the uniform strong law
fails for a stationary ergodic process. These two papers give ample
background and provide illuminating discussions; they are an
excellent starting point for anyone wishing to delve deeper into the
topic.
%interested
%reader is encouraged to consult them.

%etc... won't review all techniques here, refer to my thesis, Chatterjee's
%thesis, Ledoux's book...

%The main case of
%interest for this paper is $\R^n$ with the Lebesgue $\sigma$-algebra.
%
%
%The second goal of this paper is expository: to familiarize the
%statistical community with some recent concentration bounds for
%nonproduct measures, and how these can be applied to obtain laws of
%large numbers for interesting processes (such as Markov Chains).

This paper is organized as follows. In Section
\ref{sec:notation} we set down the notations and definitions used
throughout the paper. This is followed by Section
\ref{sec:conc-martingale}, where we describe the martingale method as
our
main workhorse for proving concentration of measure results. Our main
law of large numbers is stated in Section \ref{sec:slln}. This strong
law of large numbers is applied to adaptive Markov chains in Section 
\ref{sec:mmc}. Finally, some technical lemmas are deferred to the Appendix.
%{\bf [to be filled in later]}

%\item state LLN
%\item follows from measure concentration
%\item state bound
%\een
%\item Techniques
%\ben
%\item Martingale difference
%\item $\eta$-mixing
%\item $\Phi$, $\Psi$ norms, inequality
%\item Martingale bound
%\item concentration bound
%\een
%\item Applications
%\ben
%\item (Weakly) Contracting Markov Chains
%\item Hidden Markov Chains
%\een
%\item Technical lemmas
%\ben
%\item Minorization implies contraction
%\item $\eta$ coefficients of a Markov chain
%\item $\eta$ coefficients of a Hidden Markov chain
%\een
%\een

\section{Notation and Definitions}
\label{sec:notation}

Let
$(\X^n,\calF^n,\P)$
be a probability space, where
$\calF^n$ is the usual Borel $\sigma$-algebra generated by the finite
dimensional cylinders.
On this space
define the random process
$(X_i)_{1\leq i\leq n}$,
$X_i\in\X$.
Throughout this paper, we assume $\P\abscont\mu^n$, for some positive
Borel product measure $\mu^n=\mu\otimes\mu\otimes\ldots\otimes\mu$
on $(\X^n,\calF^n)$.

We use the indicator variable
$\pred{A}(\omega)$, equal to 1 if $\omega \in A$ and 0 otherwise,
 for $\omega \in \X^n$.   For brevity, we will occasionally omit
the argument $\omega$.
%The sign function can then be defined by
%$\sgn(z)=\pred{z>0}-\pred{z<0}$,
%We also use
%and
The ramp function
is defined
by $\pl{z}=z\pred{z>0}$. 
For arbitrary $\X$, we
%may
define the
{\em Hamming} metric on $\X^n$:
\beqn
\label{eq:ham}
d_\ham(x,y) &=& \sum_{i=1}^n \pred{x_i\neq y_i},
\qquad x,y\in\X^n.
\eeqn

We need to introduce
the notion of $\eta$-mixing, defined in
\cite{kont06-metric}\footnote{This notion of mixing is distinct
  from, and not to be confused with $\eta$-weak dependence of
%Doukhan et al.
\cite{doukhan99}.}.
We note that this type of mixing
is by no means new; it can be
traced (at least implicitly) to
%Marton's work
\cite{marton98} and is
quite explicit in \cite{samson00} and
\cite{chazottes07}.
For $1\leq i<j\leq n$
and $x\in\X^i$,
let
$$\calL(\sseq{X}{j}{n}\gn \sseq{X}{1}{i}=x)
$$ be the law
(distribution)
of $\sseq{X}{j}{n}$ conditioned on
$\sseq{X}{1}{i}=x$.
For $y\in\X^{i-1}$ and
$w,w'\in\X$,
define
\beqn
\label{eq:hdef}
\eta_{ij}(y,w,w') &=&
\TV{
\calL(\sseq{X}{j}{n}\gn \sseq{X}{1}{i}=\scat{y\sd w})-
\calL(\sseq{X}{j}{n}\gn \sseq{X}{1}{i}=\scat{y\sd w'})
},
\eeqn
where for signed measures $\nu$,
$\TV{\nu} = [\sup_A \nu(A) - \inf_A \nu(A)]/2$
is the total variation norm and
\beq
\bar\eta_{ij} &=&
\essup_{y\in\X^{i-1},w,w'\in\X}
\eta_{ij}(y,w,w'),
\eeq
the essential supremum being taken with respect to
%the marginal measure of $\P$ on $\X^i$.
$\mu$.
\hide{
Recall
that if
$(U,\cal U,\P)$,
is a probability space
and
$f:U\to\R^+$
is measurable,
$
\essup_{x\in U}f(x)$
is the smallest $a\in[0,\infty]$
for which
$f\leq a$ holds $\P$-almost surely;
this quantity is also written as $\nrm{f}_{L_\infty(\P)}$.
}

\begin{rem}
The definition in (\ref{eq:hdef}) is flawed as stated, as it does not adequately handle the case where 
sets of measure zero
are
conditioned
on. 
A more precise definition is given in \citet{kont08-pres}.
%One could attempt to perturb $\mu$ slightly so as to make it have full support and appeal to continuity
%arguments, but here too care must be taken. 
Taking
$\calP_+^n(\X)$ 
to
be the set of all 
strictly positive
probability measures $\mu$ on $\X^n$ 
(i.e.,
$\mu(x)>0$ for all $x\in\X^n$),
it 
is not hard to show 
that the functional
$\etab_{ij}: \calP_+^n(\X) \to \R$ is continuous 
on $\calP_+^n(\X)$
with respect to $\TV{\cdot}$. 
However, 
this continuity can break down
on the boundary of $\calP_+^n(\X)$. 
Thus, we may arbitrarily define $\etab_{ij}$ on any probability measure $\mu$ by
\beqn
\label{eq:hinf}
\etab_{ij}(\mu) &=& \inf_{\set{\mu_k}} \lim_{k\to\infty} \etab_{ij}(\mu_k)
\eeqn
where the infimum is taken over all sequences $\set{\mu_k : \mu_k\in \calP_+^n(\X), \TV{\mu_k-\mu}\to0}$.

See Section 5.4 of \citet{kont07-thesis}
for a discussion of the continuity of $\etab$ and conditioning on sets of measure zero, as well as motivation
for the definition in (\ref{eq:hinf}).
\end{rem}

Let $\Dn$ be the upper-triangular $n\times n$ matrix defined by
$(\Dn)_{ii}=1$ and
\beqn
\label{eq:Ddef}
(\Dn)_{ij} = \bar\eta_{ij}.
\eeqn
for $1\leq i<j\leq n$. Recall that the $\ell_\infty$
operator norm is given by
\beqn
\tsnrm{\Dn}_\infty
&=&
\label{eq:infnorm}
\max_{1\leq i< n}(
1 + \bar\eta_{i,i+1} + \ldots +\bar\eta_{i,n}
).
\eeqn

We will occasionally want to make the dependence of
$\Dn$ on the probability measure $P$ explicit; to do this
we will write $\Dn(P).$
Let us collect some simple observations about
$\nrm{\Dn(\cdot)}_\infty$:
\belen
Let $P \abscont\mu^n$ be a probability measure on
$(\X^n,\calF)$.
Then
\bit
\item[(a)] $1\leq \nrm{\Dn(P)}_\infty \leq n$
\item[(b)] $\nrm{\Dn(P)}_\infty=1$ iff $P$ is (equivalent
  a.e. $[\mu^n]$ to) a product measure
\item[(c)] if $Q\abscont\mu^m$ is a probability measure on $\X^m$ then
\beq
\nrm{\Delta_{m+n}(P\tp Q)}_\infty &\leq&
\max\set{\nrm{\Delta_{n}(P)}_\infty,\nrm{\Delta_{m}(Q)}_\infty}.
\eeq
\eit
\enlen
These properties are established in \cite{kont07-thesis},
which also gives
%As a trivial observation, note that if the variables $(X_i)$ are
%mutually independent, we have
%$(\Dn)_{ij}=\pred{i=j}$ and
%$\nrm{\Dn}_\infty=1$.
%See \cite{kont06-metric} for
a discussion of the relationship between $\eta$-mixing
and $\phi$- and other kinds of mixing.

\section{Concentration via Martingale Differences}
\label{sec:conc-martingale}

%\label{sec:mgale-bkg}
Recall our probability space $(\X^n,\calF^n,\P)$ and
let $\calF_i$ be the $\sigma$-algebra generated by
$(\seq{X}{1}{i})$,
which induces the filtration
\beqn
\label{eq:filtr}
\{\emptyset,\X^n\} = \calF_0 \subset\calF_1 \subset\ldots
\subset\calF_n = \calF^n.
\eeqn
For $i=1,\ldots,n$ and
$f\in L_1(\X^n,\P)$,
define the martingale difference
\beqn
\label{eq:mdifdef}
V_i &=& \E[f\gn \calF_i] - \E[f\gn \calF_{i-1}].
\eeqn
It is a classical result\footnote{
See \cite{ledoux01} for a modern presentation and a short proof of (\ref{eq:azuma}).
}, going back to \cite{azuma}, that
\beqn
\label{eq:azuma}
\pr{\abs{f-\E f}>t} &\leq& 2\exp\left( {-t^2 \over
2\sum_{i=1}^n\nrm{V_i}_{\infty}^2   } \right).
\eeqn
%(the meaning of $\nrm{V_i}_{L_\infty}$ will be made explicit later).
Thus, if we are able to
uniformly bound the martingale difference,
\beq
\max_{1\leq i\leq n} \nrm{V_i}_{\infty} &\leq& H_n,
\eeq
we obtain the concentration inequality
\beqn
\label{eq:genconc}
\pr{\abs{f-\E f}>t} &\leq& 2\exp
\paren{-\frac{t^2}{2 n H_n^2}}
%(-{\ts\oo2}\eps^2 n/H^2)
.
\eeqn

%\section{Obtaining L\'evy families from continuity and mixing}
\newcommand{\yi} {\sseq{y}{1}{i-1}}
\newcommand{\yii}{\sseq{y}{1}{i}}
\newcommand{\si} {\sseq{x}{1}{i-1}}
\newcommand{\sii}{\sseq{x}{1}{i}}
\newcommand{\Si} {\sseq{X}{1}{i-1}}
\newcommand{\Sii}{\sseq{X}{1}{i}}
\newcommand{\xin}{\sseq{x}{i+1}{n}}
\newcommand{\Xin}{\sseq{X}{i+1}{n}}

\hide{
\subsection{Simple bound on the martingale difference}
\label{sec:simplemartbd}
Let $(\X^n,\calF,\P)$
be a probability space
and
$(X_i)_{1\leq i\leq n}$
%be as defined in \S\ref{sec:conc-martingale}.
its associated random process; define the
filtration $\set{\calF_i}$ as in (\ref{eq:filtr}).
}

%At this point, we make the additional
Recall our
assumption that
$d\P(x)=p(x)d\mu^n(x)$ for some
%``good''
positive Borel product
measure
$\mu^n=\mu\otimes\mu\otimes\ldots\otimes\mu$
on $(\X^n,\calF)$.
%which we refer to as the {\em carrying measure}.
Similarly, the conditional probability
satisfies
$\P(\cdot\gn\calF_i)\abscont\mu^{n-i}$, with density
$p(\cdot\gn \sseq{X}{1}{i}=\sseq{x}{1}{i})$.
Here and below
$p(\sseq{x}{j}{n}\gn \sseq{x}{1}{i})$
will occasionally be used
in place of
$p(\sseq{x}{j}{n}\gn \sseq{X}{1}{i}=\sseq{x}{1}{i})$;
no
ambiguity should
arise.

For $f\in L_1(\X^n,\P)$,
$1\leq i\leq n$
and $\sseq{y}{1}{i}\in\X^i$, define
\beqn
\label{eq:videf}
V_i(f;\sseq{y}{1}{i}) &=&
\E[f(X)\gn\sseq{X}{1}{i}=\sseq{y}{1}{i}] -
\E[f(X)\gn\sseq{X}{1}{i-1}=\sseq{y}{1}{i-1}];
\eeqn
this is just the martingale difference.

A slightly more tractable
quantity turns out to be
\beqn
\label{eq:vhatdef}
\hat V_i(f;\yi,w_i,w_i') &=&
\E[f(X)\gn\Sii = \scat{\yi \sd w_i}]-
\E[f(X)\gn\Sii = \scat{\yi \sd w_i'}],
\eeqn
where
$w_i,w_i'\in\X$.
These two quantities have a simple relationship, which may be stated
symbolically as
\beqn
\label{eq:vhat}
\tsnrm{V_i(f;\cdot)}_{L_\infty(\P)} &\leq &
\tsnrm{\hat V_i(f;\cdot)}_{L_\infty(\P)};
\eeqn
this is proved in \cite[Lemma 4.1]{kont06-metric}.

The next step is to notice that
$\hat V_i(\cdot;\yi,w_i,w_i')$, as a functional on $L_1(\X^n,\P)$, is
linear; in fact, it is given by
\beqn
\label{eq:linfunc}
\hat V_i(f;\yi,w_i,w_i') \;=\;
\int_{\X^n} f(x) \hat g(x)d\mu^n(x)
\;\defeq\;
\iprod{f,\hat g},
\eeqn
where
\beqn
\label{eq:ghat0}
\hat g(x) &=&
\pred{\sii= \scat{\yi\sd w_i}}
p(\xin\gn \scat{\yi\sd w_i})
-
\pred{\sii =\scat{\yi\sd w_i'}}
p(\xin\gn \scat{\yi\sd w_i'}).
\eeqn
The plan is to bound $\iprod{f,\hat g}$ using continuity properties
of $f$ and mixing properties of $X$, which will immediately lead to a
result of type (\ref{eq:genconc}) via (\ref{eq:vhat}).

\subsection{$\Phi$ and $\Psi$ Norms}
To state our results in sufficient generality, we shall borrow several
definitions from \cite{kont07-thesis}.
Let $(\calX,\rho)$ be a metric space and
recall the definition of the {\em Lipschitz constant} of an
$f:\calX\to\R$:
\beqn
\Lip{f} =
\sup_{x\neq y} \frac{\abs{f(x)-f(y)}}{\rho(x,y)},
\qquad x,y\in\calX.
\eeqn
%Note that Lipschitz-continuity does not guarantee measurability
%(indeed, no measure has been specified as of yet).
Recall also the
definition of the {\em diameter}:
\beq
\diam_{\rho}(\calX) &=& \sup_{x,y\in\calX}\rho(x,y).
\eeq

Let $\mu$ be a positive Borel measure on
a measurable space
$(\X,\calF)$ and
let $F_n=L_1(\X^n,\calF^n,\mu^n)$
be equipped
with the inner product
\beqn
\label{eq:ipdef}
 \iprod{f,g} &=& \int_{\X^n} f(x)g(x) d\mu^n(x).
\eeqn
Since $f,g\in F_n$ might not be in $L_2(\X^n,\calF^n,\mu^n)$, the expression in
(\ref{eq:ipdef})
in general
might not be finite. However, for
%$g\in\Lips(\X^n,\rho)$,
$g\in L_\infty(\X^n,\calF^n,\mu^n)$,
we
have
\beqn
\label{eq:phitriv}
\abs{\iprod{f,g}} &\leq&
%\diam_{\rho}(\X^n)
\nrm{f}_{1}
\nrm{g}_\infty
.
\eeqn

Define
%The continuous analog of
%As above, define
the {\em projection} operator $\pi:F_{n}\to
F_{n-1}$
%is defined
as follows. If
$f:\X^n\to\R$ then
$(\pi f):\X^{n-1}\to\R$ is given by
\beqn
\label{eq:pidef}
(\pi f)
(x_2,\ldots,x_n) &\defeq& \int_\X f(x_1,x_2,\ldots,x_n) d\mu(x_1).
\eeqn
Note that by Fubini's theorem (Thm. 8.8(c) in \cite{rudin87}),
$\pi f \in L_1(\X^{n-1},\mu^{n-1})$.
Define the functional $\psin{n}:F_n\to\R$
recursively: $\psin{0}\defeq 0$ and
\beqn
\label{eq:contpsidef}
 \psin{n}(f) &\defeq&
\psin{n-1}(\pi f)
+
\int_{\X^n} \pl{f(x)}d\mu^n(x)
\eeqn
for $n\geq1$.
The latter is finite since
\beqn
\label{eq:psitriv}
\psin{n}(f) &\leq& n\nrm{f}_{L_1(\mu)},
\eeqn
\hide{
Clearly $(\X^n,\sigalg,\mu^n)$ is a measure space.
Let $F_n=L_1(\X^n,\mu^n)$ be the usual space
of integrable functions $f:\X^n\to\R$
and define $\Phi_n\subset F_n$:
\beq
\Phi_n &=& \set{f:\X^n\to[0,n] ; \Lip{f}\leq 1}
\eeq
where
\beqn
\Lip{f} =
\sup_{x\neq y} \frac{\abs{f(x)-f(y)}}{d(x,y)},
\qquad x,y\in\X^n,
%\inf \{ c \in \R :
%\abs{f(x)-f(y)} &\leq& c d_\ham(x,y)~\forall~x,y\in \X^n \}
\eeqn
with
\beqn
\label{eq:ham}
d(x,y) &=& \sum_{i=1}^n \pred{x_i\neq y_i},
\eeqn
denoting the Hamming metric.
Note that $\Lip{f}\leq 1$ does not guarantee that $f$ is
$\sigalg$-measurable,
so the requirement that $\Phi_n\subset F_n$ is essential.

Equip $F_n$
with the inner product
\beqn
\label{eq:ipdef}
 \iprod{f,g} &=& \int_{\X^n} f(x)g(x) d\mu^n(x).
\eeqn
Since $f,g\in F_n$ might not be in $L_2(\X^n,\mu^n)$, the expression in
(\ref{eq:ipdef})
in general
%might
will
not be finite. However, for $g\in\Phi_n$, we
have, by H\"older's inequality,
\beqn
\label{eq:phitriv}
\abs{\iprod{f,g}} &\leq& n \nrm{f}_{L_1}
\eeqn
(the motivation for
bounding
$\iprod{f,g}$
comes from
(\ref{eq:linfunc})).

Define the {\em marginal projection} operator $\pi:F_{n}\to F_{n-1}$ as
follows. If
$f:\X^n\to\R$ then
$(\pi f):\X^{n-1}\to\R$ is given by
\beqn
\label{eq:pidef}
%h
(\pi f)
(x_2,\ldots,x_n) &\defeq& \int_\X f(x_1,x_2,\ldots,x_n) d\mu(x_1).
\eeqn
Note that by Fubini's theorem (e.g. Thm. 8.8(c) in \cite{rudin87}),
%$\pi f$ is $\mu$-measurable.
$\pi f \in L_1(\X^{n-1},\mu^{n-1})$.
Define the functional $\psin{n}:F_n\to\R$
recursively: $\psin{0}\defeq 0$ and
\beqn
\label{eq:contpsidef}
 \psin{n}(f) &\defeq& \psin{n-1}(\pi f) +\int_{\X^n} \pl{f(x)}d\mu^n(x)
\eeqn
for $n\geq1$.
%Observe that
The latter is finite since
\beqn
\label{eq:psitriv}
\psin{n}(f) &\leq& n\nrm{f}_{L_1},
\eeqn
}
as shown in the following Lemma.
%\cite[Thm. A.1]{kont06-metric}.

\hide{
We are about to define two functionals on
$F_n$.
Although we use the norm
notation, none of the results we prove actually rely on the norm
properties of
$\phinorm{\cdot}$ and
$\psinorm{\cdot}$.
It turns out
that under appropriate conditions both
are valid norms \cite{kont06-metric};
$\psinorm{\cdot}$ is (topologically) equivalent to $\nrm{\cdot}_{L_1}$ while
$\phinorm{\cdot}$ is weaker. NO!!
}

Let $\Phi_n\subset F_n$ be the set of all
measurable\footnote{
Note that $\Lip{f}\leq 1$ does not guarantee that $f$ is
$\sigalg$-measurable,
so the requirement that $\Phi_n\subset F_n$ is essential.}
$f:\X^n\to[0,\diam_\rho(\X^n)]$ with $\Lip{f}\leq1$ and
define two norms on $F_n$:
\beqn
\label{eq:phinorm}
\phinorm{f} &=& \sup_{g\in \Phi_n} \abs{\iprod{f,g}}
\eeqn
and
\beqn
\label{eq:psinorm}
 \psinorm{f} &=& \max_{s\in\set{-1,1}}\psin{n}(s f).
\eeqn
We refer to the norms in
(\ref{eq:phinorm})  and (\ref{eq:psinorm})
as $\Phi$-norm and $\Psi$-norm, respectively, and summarize some of
their properties:
\belen
\label{lem:phipsiprops}
Under mild measure-theoretic regularity conditions on
$(\X^n,\sigalg,\mu^n)$, which cover the case of
%$\rho=d_\ham$ and
$\X$ countable and
$\X=\R$ with Lebesgue measure
(the metric being
$\rho=d_\ham$ in either case),
the functionals defined in
(\ref{eq:phinorm})
and
(\ref{eq:psinorm})
satisfy
\bit
\item[(a)]
$\Phi$-norm and $\Psi$-norm are valid vector-space norms on $F_n$
\item[(b)]
for all $f\in F_n$,
\beq
          {\ts\oo2}\nrm{f}_{L_1} \;\leq \; \phinorm{f}
\;\leq \; n        \nrm{f}_{L_1}
\eeq
\item[(c)]
for all $f\in F_n$,
\beq
          {\ts\oo2}\nrm{f}_{L_1} \;\leq \; \psinorm{f}
\;\leq \; n        \nrm{f}_{L_1}
\eeq
\eit
\enlen
\bepf
%\bit
%\item[(a)] This
The claim in
(a) is proved in
\cite[Theorems A.1, A.2, A.3]{kont06-metric};
(b) and (c) are proved ibid. in (52) and Theorem A.1(b), respectively.
%\item[(b)]
%\eit
\enpf

\begin{definition}
A metric space $(\X^n,\rho)$ is
said to be
{\em $\Psi$-dominated}
with respect to
a positive Borel measure $\mu$ on $\X$
if
the inequality
\beqn
\label{eq:psidom}
\sup_{g\in\Phi_n} \iprod{f,g} &\leq& \psin{n}(f)
\eeqn
holds
for all $f\in F_n$.
\end{definition}

\begin{rem}
We have defined $\Phi_n$, $\Psi_n(\cdot)$, and $\iprod{\cdot,\cdot}$ in an
abstract measure space. To indicate explicitly what space they are
being defined over, we will use the notation
$\phix{\X^n,\mu^n}$, $\psix{\X^n,\mu^n}(\cdot)$, and
$\iprodx{\X^n,\mu^n}{\cdot,\cdot}$.
\end{rem}

When $\X$ is a finite set with counting measure $\nu$,
and $\rho=d_\ham$,
we have
\beqn
\label{eq:phipsidiscr}
\iprodx{\X^n,\nu^n}{f,g} &\leq & \psix{\X^n,\nu^n}(f)
\eeqn
for all $f:\X^n\to\R$ and $g\in\phix{\X^n,\nu^n}$;
in other words, $(\X^n,d_\ham)$ is $\Psi$-dominated with respect to $\nu$.
This was first
proved in \cite{kontram06}; see \cite{kont06-lp} for a much
simpler proof.
\cite{kont07-thesis} extends this to countable $\X$.
The goal of the remainder of this section is to establish the analogue of
(\ref{eq:phipsidiscr}) for $\X=\R$ with the Lebesgue measure $\mu$.

\bethn
\label{thm:phipsi}
Let $\mu^n$ be the Lebesgue measure on $\R^n$
and take $\rho=d_\ham$. Then $(\R^n,d_\ham)$ is $\Psi$-dominated.
That is, for all
$f\in L_1(\R^n,\mu^n)$ and $g\in\phix{\R^n,\mu^n}$, we have
\beqn
\label{eq:phipsiR}
\iprodx{\R^n,\mu^n}{f,g} &\leq & \psix{\R^n,\mu^n}(f).
\eeqn
\begin{rem}
Observe that (\ref{eq:phipsiR}) is equivalent to
\beqn
\label{eq:phipsiRnrm}
\phinorm{f} &\leq&\psinorm{f},
\qquad f\in L_1(\R^n,\mu^n).
\eeqn
The proof will closely follow the argument in
\citet[Theorem 8.1]{kont06-metric}.
\end{rem}
\bepf
Let $C_c$ denote the space of continuous functions $f:\R^n\to\R$ with
compact support;
it follows from
\cite[Theorem 3.14]{rudin87} that $C_c$ is dense in $\Lone$,
in the topology induced by $\nrm{\cdot}_{L_1}$.
This
implies
that
for any $f\in\Lone$
and $\eps>0$, there is a $g\in C_c$ such that
$\nrm{f-g}_{L_1}<\eps/n$
and therefore
(via Lemma
\ref{lem:phipsiprops}(b) and (c)),
\beq
\phinorm{f-g}<\eps
\qquad\text{and}\qquad
\psinorm{f-g}
< \eps
\eeq
so it suffices to prove
(\ref{eq:phipsiRnrm}) for $f\in C_c$.

For $m\in\N$, define $Q_m\subset \QQ$ to be the rational numbers with denominator
$m$:
\beq
Q_m &=& \set{p/r\in\QQ : r=m}.
\eeq
Define the map $\gamma_m:\R\to Q_m$ by
\beq
\gamma_m(x) &=& \max\set{q\in Q_m : q\leq x}
\eeq
and extend it to $\gamma_m:\R^n\to Q_m^n$ by defining
$[\gamma_m(x)]_i = \gamma_m(x_i)$. The set $Q_m^n\subset\R^n$ will be
referred to as the $m$-{\em grid points}.

We say that $g\in\Lone$ is a
{\em grid-constant function} if there is an $m>1$ such that $g(x)=g(y)$
whenever $\gamma_m(x)=\gamma_m(y)$; thus a grid-constant
function is constant on the grid cells induced by $Q_m$.
%for some grid.
Let $G_c$ be the space of the grid-constant functions with compact
support; note that $G_c\subset\Lone$.
It is easy to see that $G_c$ is
dense in $C_c$.
Indeed,
%for
pick any
$f\in C_c$
and let $M\in\N$ be such that
$\supp(f)\subset[-M,M]^n$.
Now a continuous function is uniformly continuous on a compact set,
and so for any
$\eps>0$, there is a $\delta>0$
such that $\omega_f(\delta)<\eps/(2M)^n$, where $\omega_f$ is the
$\ell_\infty$
modulus of
continuity of $f$.
Take $m=\ceil{1/\delta}$
%and $M\in Q_m$ to be such that $\supp(f)\subset[-M,M]^n$;
and let $g\in G_c$ be
such that
$\supp(g)\subset[-M,M]^n$ and
$g$ agrees with $f$ on the $m$-grid points.
Then
we have
\beq
\nrm{f-g}_{L_1} \leq
(2M)^n\nrm{f-g}_{L_\infty} < \eps.
\eeq
Thus we need only prove
(\ref{eq:phipsiR})
for $f\in G_c$, $g\in G_c\cap \phix{\R^n,\mu^n}$.

Let
$f\in G_c$ and $g\in G_c\cap \phix{\R^n,\mu^n}$ be fixed, and let $m>1$
be such that $f$ and $g$ are $m$-grid-constant functions.
Let
$\bar\k,\bar\f:Q_m^n\to\R$
be such that
$\bar\k(\gamma_m(x))=f(x)$
and
$\bar\f(\gamma_m(x))=g(x)$
for all $x\in\R^n$. As above, choose $M\in\N$ so that
$\supp(f)\cup\supp(g)\subset[-M,M]^n$. Then, denoting the counting
measure on $Q_m^n$ by $\nu^n$, we have
\beq
\iprodx{\R^n,\mu^n}{f,g} &=&
\paren{2M \over m}^n \iprodx{Q_m^n,\nu^n}{\bar\k,\bar\f}
\eeq
and
\beq
\psix{\R^n,\mu^n}(f) &=&
\paren{2M \over m}^n \psix{Q_m^n,\nu^n}(\bar\k).
\eeq
Now $Q_m$ is finite and by construction,
$\bar\f\in\phix{Q_m^n,\nu^n}$, so
(\ref{eq:phipsidiscr}) applies.
This shows
$\iprodx{\R^n,\mu^n}{f,g} \leq
\psix{\R^n,\mu^n}(f)$ and completes the proof.
\enpf

\enthn

\subsection{Bounding the martingale difference}
The machinery of
$\eta$-mixing
and
%the ``$\Phi-\Psi$'' inequality in Theorem \ref{thm:phipsi}
$\Psi$-dominance
allows us to bound the martingale difference for a Lipschitz function
of arbitrarily dependent random variables.

\bethn
\label{thm:mgalebd}
Let $(\X,\calF,\mu)$ be a positive Borel measure space and suppose
that $(\X^n,\rho)$ is $\Psi$-dominated with respect to $\mu$ for some
metric $\rho$ on $\X^n$. Let $(\X^n,\calF^n,\P)$ be a\hide{
Let $(\R^n,\sigalg,\mu^n)$ be the usual $n$-dimensional Lebesgue
measure space
and
%let
suppose
$(\R^n,\sigalg,\P)$
%be
is
a}
probability space
with
$\P\abscont\mu^n$.
Then for $1\leq i\leq n$,
\beqn
\label{eq:Vhbound}
\nrm{V_i(f;\cdot)}_{L_\infty(\P)} &\leq&
\Lip{f}\nrm{\Dn(\P)}_\infty,
\eeqn
where $V_i$ is the martingale difference defined in
(\ref{eq:mdifdef})
%, the Lipschitz constant $\Lip{f}$ is defined with
%respect to the Hamming metric (\ref{eq:ham}),
and
$\Dn$ is the $\eta$-mixing matrix defined in
(\ref{eq:Ddef}).
\enthn
\begin{rem}
This result can be proved, almost verbatim, by the argument given in
\citet[Theorem 7.1]{kont06-metric}, so we only give a sketch of the proof here.
\end{rem}
\bepf
Since
$\tsnrm{V_i(f;\cdot)}_{L_\infty}$
and
$\Lip{f}$
are both homogeneous functionals of $f$ (in the sense of
$T(af)=|a|T(f)$ for $a\in\R$), there is no loss of generality in
taking $\Lip{f}=1$.
Additionally, since
$V_i(f;y)$ is
translation-invariant (in the sense that $V_i(f;y)=V_i(f+a;y)$ for
all $a\in\R$), there is no loss of generality in restricting the
range of $f$ to $[0,n]$. In other words, it suffices
to consider $f\in\phix{\X^n,\mu^n}$.

From (\ref{eq:vhat}) we have that it suffices to bound
$\tsnrm{\hat V_i(f;\cdot)}_{L_\infty}$,
defined in (\ref{eq:vhatdef}).
By (\ref{eq:linfunc}),
%having an
we have the
equivalent form
\beqn
\label{eq:mgaleiprod}
\hat V_i(f;\yi,w_i,w_i') \;=\;
\int_{\X^n} f(x) \hat g(x)d\mu^n(x)
\;=\;
\iprodx{\X^n,\mu^n}{f,g_i},
\eeqn
where
$g_i$
has a simple explicit construction (\ref{eq:ghat0}),
depending on $\yi,w_i,w_i'$.

It is shown
in the course of proving
\cite[Theorem 7.1]{kont06-metric}
that
\beq
\iprod{f,g_i} &=& \iprod{T_y f,T_y g_i},
\eeq
where
the operator
$T_y:\Lone\to {L_1(\X^{n-i+1},\mu^{n-i+1})}$
is
defined by
\beq
(T_y h)(x) &\defeq& h(yx),
\qquad\text{for all } x\in\X^{n-i+1}.
\eeq

Appealing to
Theorem \ref{thm:phipsi},
we get
\beqn
\iprod{T_y f,T_y g_i} &\leq & \psin{n}(T_y g_i)
.
\eeqn

Furthermore,
as shown ibid.,
the form of $g_i$ implies that
\beqn
\psin{n}(T_y g_i) &\leq& 1+ \sum_{j=i+1}^n \bar\eta_{ij},
\eeqn
establishing (\ref{eq:Vhbound}).

\enpf

\section{The Strong Law of Large Numbers}
\label{sec:slln}

We are now in a position to state our main result.

\bethn
\label{thm:lln}
Let $(\X,\calF,\mu)$ be a positive Borel measure space and suppose
that $(\X^n,d_\ham)$ is $\Psi$-dominated with respect to $\mu$.
% for some
%metric $\rho$ on $\X^n$.
Define
the random process
$\sseq{X}{1}{\infty}$
%be a random process defined
on the measure
space
$(\X^\N,\sigalg^\N,\P)$, and assume that
for all $n\geq1$ we have
$\Pn\abscont\mu^n$,
where $\Pn$ is
the marginal distribution on $\sseq{X}{1}{n}$ and $\mu^n$ is
the corresponding product measure on
$(\X^n,\sigalg^n)$.
%the
%Lebesgue measure on $(\R^n,\sigalg^n)$.
Suppose further that
the empirical
measure defined by
\beqn
\label{eq:pndef}
\hat P_n(A) &=& \oo n \sum_{i=1}^n \pred{X_i \in A},
\quad A \in \sigalg,
\eeqn
has
uniformly
converging expectation:
\beq
%\lim_{n\to\infty} \mexp \hat P_n(A)
%&=& \nu(A);
\lim_{n\to\infty} \TV{\mexp \hat P_n(\cdot)
- \nu(\cdot)}
&\to& 0
\eeq
and
define $n_0=n_0(\eps)$ to be such that
%$\abs{\mexp \hat P_n(A)-\nu(A)}<\eps$
$\TV{\mexp \hat P_n(\cdot)-\nu(\cdot)}<\eps$
for all $n>n_0(\eps)$.

Then $\hat P_n(A)$ converges to $\nu(A)$ almost surely, exponentially
fast:
\beqn
\label{eq:expfast}
%\pr
\Pn\tlprn
{\abs{\hat P_n(A)-\nu(A)}>t+\eps} &\leq&
2\exp(-nt^2/2\nrm{\Dn}_\infty^2)
\eeqn
for all $n>n_0(\eps)$,
where
$\Dn$ is the
$\eta$-mixing matrix defined in equation~(\ref{eq:Ddef}).
\enthn

Particular cases of interest include $\X$ countable with
$\mu$ taken to be counting measure and
$\X=\R$ with $\mu$ taken to be Lebesgue measure.

\bepf
This result follows directly from Theorem \ref{thm:mgalebd}, by
observing that
the function
$\f_A:\sseq{X}{1}{n}\to\R$ defined by
$\f_A(\sseq{X}{1}{n})=\hat P_n(A)$ has Lipschitz constant $1/n$.
\enpf

\becon
\label{cor:asconv}
Under the conditions of Theorem \ref{thm:lln}, $\hat P_n$ converges to $\nu$ in 
%$\TV{\cdot}$, 
%total variation,
distribution,
almost surely.
\encon
\bepf
This is an
immediate consequence of the (first) Borel-Cantelli lemma.
\enpf

\section{Concentration of Marginals of Markov chains}
\label{sec:mmc}
For clarity of presentation, we take all the state spaces to be finite, until indicated otherwise in
Section \ref{sec:apps}.
%For clarity, we present the calculations in this section over finite sets; 
Everything extends
easily to the continuous case, as shown in the Appendix.

\subsection{Bounding $\eta_{ij}$ via contraction coefficients}
\label{sec:eta-theta}

Consider a Markov process
$$
(W_t)_{t=1,2,\ldots}
$$
taking values $W_t = (X_t, Y_t) \in \observedX \times \hiddenX,$
defined on the probability space
$( (\observedX \times \hiddenX)^\N,
(\sigalg_{\osymb} \tp \sigalg_{\hsymb})^\N , \P )$;
the subscripts `o' and `h' are used to suggest ``observed'' and ``hidden'' states.
Suppose that we are interested primarily in the marginal
behaviour of $(X_t)$.  (This might be the case, for instance,
if we can observe $(X_t)$ but not $(Y_t)$ as when analyzing
hidden Markov models.  It is also of interest in the
context of analysis of adaptive Markov chain Monte Carlo schemes
as described later in this section.)
Let us write the transition kernel of $(W_t)$ as
$$
K_t(w,A) = \P(W_{t+1} \in A \gn W_t=w),
\quad
w \in (\observedX \times \hiddenX),
A \in (\sigalg_{\osymb} \tp \sigalg_{\hsymb}),
$$
and assume that the initial distribution of the process is
$$
p_0(A) = \P(W_1 \in A).
$$

When $(X_t)_{t \in \N}$
can be embedded into a higher-dimensional
Markov process as stated above, we will
call $(X_t)$ a \emph{Markov marginal chain} (MMC).
Note that Markov marginal chains properly contain
the ordinary 
%and hidden 
Markov chains, and are easily seen to be equivalent in expressive
power to hidden Markov chains.

In this section we apply the results established
in the previous sections to obtain relatively simple conditions
for convergence of empirical measures associated
with the MMC.

\hide{
Consider two finite sets, $\hiddenX$
(the ``hidden'' state-space)
and $\observedX$
(the ``observed'' state-space).

Let $\mu$ be a Markov measure
on $(\observedX\times\hiddenX)^n$
defined by the initial distribution $p_0$ and the kernels
$\set{K_i(\cdot\gn\cdot)}_{1\leq i\leq n}$:
\beqn
\label{eq:mudef}
\mu\paren{\ohs{\sseq{\observedx}{1}{n}}{\sseq{\hiddenx}{1}{n}}}
&=&
p_0\paren{\oh{x_1}{x_1}}
\prod_{i=1}^{n-1}
K_i\paren{\oh{x_{i+1}}{x_{i+1}}
\Given \oh{x_{i}}{x_{i}}
},
\eeqn
where for readability we use the stacked notation $\paren{\oh{x}{x}}$
--
instead
of the more standard $(\observedx,\hiddenx)$
--
for elements of $\observedX\times\hiddenX$.
}

\hide{
Recall the
%definition of
$\eta$-mixing coefficients, defined on some
measure space $(\calX^n,\P)$.
%$\P$ over $\calX^n$.
For $1\leq i<j\leq n$
and $x\in\calX^i$,
let
$$\calL(\sseq{X}{j}{n}\gn \sseq{X}{1}{i}=x)
$$ be the law
(distribution)
of $\sseq{X}{j}{n}$ conditioned on
$\sseq{X}{1}{i}=x$.
For $y\in\calX^{i-1}$ and
$w,w'\in\calX$,
define
\beq
\eta_{ij}(y,w,w') &=&
\TV{
\calL(\sseq{X}{j}{n}\gn \sseq{X}{1}{i}=y w)-
\calL(\sseq{X}{j}{n}\gn \sseq{X}{1}{i}=y w')
},
\eeq
where
$\TV{\cdot}$ is the total variation norm
and
\beq
\bar\eta_{ij} &=&
\max_{y\in\calX^{i-1}}
\max_{w,w'\in\calX}
\eta_{ij}(y,w,w').
\eeq
}

Let us define the $i$\th {\em contraction coefficient} $\tha_i$ of
the MMC defined above by
\beqn
\label{eq:thadef}
\tha_i &=&
\sup_{\observedx,\observedx'\in\observedX}
\sup_{\hiddenx,\hiddenx'\in\hiddenX}
\TV{
K_i\paren{\oh{x}{x},\cdot} -
K_i\paren{\ohs{\op x}{\hp x},\cdot}
}.
\eeqn

We 
obtain a bound on
%claim that 
the 
$\eta$-mixing coefficients
of the MMC
in terms of
%can be controlled by 
its contraction coefficients:
\bethn
\label{thm:pomc}
The MMC $(X_t)$
on $\observedX^\N$, as defined above, satisfies
\beqn
\label{eq:pomc}
\bar\eta_{ij} &=&
\tha_i \tha_{i+1} \ldots \tha_{j-1}.
\eeqn
\enthn

\bepf
The proof is greatly simplified by an observation of \cite{marton07} -- namely, that a Marginal Markov chain is a special case of a Hidden Markov
chain (see \citet{rabiner89}). This is easily seen by considering the function
\beq
\pi : \paren{\oh{x}{x}} \mapsto \observed{x},
\eeq
which projects an (observed,hidden) pair onto its observed component, mapping a Markov chain to its hidden Markov marginal. It has already been shown
(see \citet{kont06-hmm} or \citet{kontram06}) that the $\eta$-mixing coefficients of a hidden Markov chain are controlled by the contraction coefficients
of the underlying Markov chain, in the manner of (\ref{eq:pomc}).
\enpf

\subsection{Tensorization lemma}
\label{sec:tensor}

\belen
\label{lem:tens}
Suppose we have a MMC on $(\observedX\times\hiddenX)^n$ defined by transition
kernels $\set{K_i(\cdot\gn\cdot)}_{1\leq i\leq n}$, which have the
following special structure:
\beqn
\label{eq:kdecomp}
K_{i}\paren{\oh{u}{u} \Given \oh{v}{v}}
&=&
A_i(\o u\gn \o v,\h v)
B_i(\h u\gn \o v,\h v).
\eeqn
Then we have
\beqn
\label{eq:thab}
\tha_i &\leq& \alpha_i+\beta_i-\alpha_i\beta_i
\eeqn
where $\tha_i$ is defined in (\ref{eq:thadef}), and
\beq
\alpha_i &=&
\max_{\observedx,\observedx'\in\observedX}
\max_{\hiddenx,\hiddenx'\in\hiddenX}
\TV{
A_i(\cdot \gn \observedx,\hiddenx)-
A_i(\cdot \gn \op x,\hp x)
},\\
\beta_i &=&
\max_{\observedx,\observedx'\in\observedX}
\max_{\hiddenx,\hiddenx'\in\hiddenX}
\TV{
B_i(\cdot \gn \observedx,\hiddenx)-
B_i(\cdot \gn \op x,\hp x)
}.
\eeq
\enlen
\bepf
The claim follows immediately from the 
well-known
total variation tensorization
lemma (see for instance
\citet{kont07-thesis}),
which states that
if $\mu,\mu'$ are probability measures on $\calX$
and $\nu,\nu'$ are probability measures on $\calY$, then
\beq
\TV{\mu\tp\nu-\mu'\tp\nu'}
&\leq&
\TV{\mu-\mu'}+\TV{\nu-\nu'}-\TV{\mu-\mu'}\TV{\nu-\nu'}
\eeq
where $\mu\tp\nu$ is a product measure on $\calX\times\calY$.
\enpf

\subsection{Concentration of adaptive Markov chains}
\label{sec:adapt}

%For clarity, we present the calculations in this section over finite sets; everything extends
%easily to the continuous case, as shown in the Appendix.
%Everything in this calculation works for continuous sets and
%integrals; I'll use sums over the finite set $\X$ for clarity.
%%Sums will range over the entire space of the summation variable;
%%thus
%%$\ds\sum_{\sseq{x}{i}{j}}f(\sseq{x}{i}{j})$ stands for
%%$$\ds\sum_{\sseq{x}{i}{j}\in\X^{j-i+1}}f(\sseq{x}{i}{j}).$$
Throughout our calculations, $n$ will be a fixed positive integer.
Let $\Gamma$ be an index set
and $\set{K_\g(\cdot\gn\cdot)}_{\g\in\Gamma}$ be a collection of
Markov transition kernels, $K_\g:\X\to\X$. For a given sequence
$\sseq{\g}{1}{n-1}\in\Gamma^{n-1}$, we have a Markov measure on
$\X^n$:
\beq
P_{\sseq{\g}{1}{n-1}}(x) &=&
p_0(x_1)\prod_{i=1}^{n-1}K_{\g_i}(x_{i+1}\gn x_i),
\qquad
x\in\X^n.
\eeq
For $1\leq i<n$, $x\in\X$, and $\g'\in\Gamma$, let $g_i(\cdot\gn x,\g')$ be a probability
measure on $\Gamma$.
Together,
$\set{K_\gamma}$ and $\set{g_i}$ define a measure $\mu$ on $\X^n$, which we call
an {\em adaptive Markov} measure:
\beqn
\label{eq:adapt-def}
\mu(x) &=& \sum_{\sseq{\g}{1}{n}\in\Gamma^{n}}
p_0(x_1,\g_1)\prod_{i=1}^{n-1}[g_i(\g_{i+1}\gn x_i,\g_{i})K_{\g_i}(x_{i+1}\gn x_{i+1})],
\eeqn
where $\g_1$ is a dummy index to make $g_1(\cdot\gn\cdot)$ well-defined.

Define the following contraction coefficients

\beqn
\label{eq:kapdef}
\kappa
&=& \max_{w,w'\in\X,\g,\g'\in\Gamma}\TV{K_{\g}(\cdot\gn w)-K_{\g'}(\cdot\gn w')}.
\eeqn

and

\beqn
\label{eq:lamdef}
\lambda_i
&=&
%\max_{1\leq i<n}
\max_{w,w'\in\X,\g,\g'\in\Gamma}
\TV{g_i(\cdot\gn w,\g)-g_i(\cdot\gn w',\g')}.
\eeqn

Then
%\paragraph{Conjecture}
\bethn
\label{thm:mmconc}
For an adaptive Markov measure $\mu$ on $\X^n$
we have
\beq
\bar\eta_{ij} &\leq& \tha_i \tha_{i+1}\cdots \tha_{j-1}
\eeq
where
%$\tha=\max\set{\kappa,\lambda}$.
$$ \tha_i = \kappa+\lambda_i - \kappa\lambda_i.$$
\enthn
\bepf
First we observe that the
adaptive Markov process
defined in (\ref{eq:adapt-def})
is in fact an MMC,
with $\observedX=\X$ and $\hiddenX=\Gamma$;
thus Theorem \ref{thm:pomc} applies.
Furthermore,
%its
the MMC
kernel
(denoted here by $P_i$ instead of $K_i$ to avoid confusion with $K_\g$)
decomposes as in (\ref{eq:kdecomp}):
\beqn
\label{eq:pg}
P_{i}\paren{\ohs{x}{\g} \Given \ohs{x'}{\g'}}
&=&
g_i(\g \gn x',\g')K_{\g'}(x \gn x').
\eeqn
The claim is proved by applying Lemma \ref{lem:tens} to (\ref{eq:pg}).
\enpf

\subsection{Example: Application to Adaptive MCMC Analysis}
\label{sec:apps}

As one application, we consider the analysis of a family of
adaptive Markov chain Monte Carlo schemes.
Such schemes have been considered in some detail by a number
of authors, beginning with studies of a specific scheme
\citep{Haarioetal, AndrieuRobert} and later being
generalized by \cite{AtchadeRosenthal}, \cite{AndrieuMoulines},
and \cite{RobertsRosenthal}.
Many other approaches to adaptation have been developed; these include, for example
\cite{599374} and \cite{citeulike:206883}.
We consider the general framework of \cite{RobertsRosenthal},
and demonstrate that by imposing a stronger form of their
so-called ``diminishing adaptation'' condition, one is able
to strengthen the weak law of large numbers they establish to a
strong law of large numbers.

Consider a stochastic process $\{X_t \in \R,~t=0,1,\ldots\}$
on $(\R^\N, \sigalg^\N, \P).$
Adopting similar notation to that of \cite{RobertsRosenthal},
define a family $\{K_\gamma(\cdot,\cdot),~\gamma \in {\cal G} \subseteq {\R}\}$
of transition kernels
such that for each $\gamma \in {\cal G}$,
$K_\gamma(\cdot,\cdot)$ is irreducible, aperiodic, and ergodic with a
limiting distribution $\pi$ on $(\R,{\cal F}).$
One would typically take each such kernel to be a Metropolis-Hastings
kernel with certain parameter values determined by $\gamma.$

For fixed $\gamma$, a homogeneous Markov chain whose joint distributions
are determined by an initial value and 
transition kernel $K_\gamma$ would have marginal distributions
converging in total variation norm to $\pi$.
However, in adaptive MCMC problems, interest centers on the behaviour
of a more complex process $\{X_t\}$.  
Rather than holding $\gamma$ fixed, one
allows the transition kernel to vary over time.
To be precise, we specify
an initial value $X_0=x_0$, along with transition probabilities
$$
P(X_{t+1} \in A \gn X_t=x) = K_{\Gamma_t}(x,A),
$$
where each 
$\Gamma_t \in {\cal G}$ is some function of $\Gamma_0,\ldots,\Gamma_{t-1},X_0,\ldots,X_t.$
Thus the kernel used at time $t$ is itself random, depending on the
past history of the process.
This means that $\{X_t\}$ is not Markovian.

To see how we can apply Theorem~\ref{thm:lln} in this context, we introduce some
assumptions.

\begin{assumption}
\label{assump:1}
There exists a random limiting kernel indexed by
$\Gamma_\infty(\omega)$ such that
\be
\Gamma_t(\omega) \rightarrow \Gamma_\infty(\omega)
\quad \forall \omega \in \Omega.
\ee
Furthermore, convergence to the random limit is uniform in
the sense that there exists a non-negative monotone
non-increasing sequence $\{\kappa_t\} \rightarrow 0$ 
such that
\be
\abs{\Gamma_t(\omega) - \Gamma_\infty(\omega)} \le \kappa_t
\quad \forall \omega \in \Omega.
\ee
\end{assumption}

\begin{assumption}
\label{assump:4}
For each $\gamma\in{\cal G}$, we have 
$K_\gamma(x,\cdot)
\convd
K_\gamma(y,\cdot)
$ 
(i.e., converges in distribution)
whenever $x\to y$.
\hide{
Each joint distribution
$\P_n$ of the adaptive chain
$\{X_1,\ldots,X_n\}$, given $X_0=x_0$,  
has a density with respect to a product
measure $\lambda^n.$
}
\end{assumption}

\begin{assumption}
\label{assump:2}
Each kernel $K_\gamma(\cdot,\cdot), \gamma \in {\cal G}$ is
uniformly ergodic, with
\be
\lim_{t \rightarrow \infty} \TV{ K_\gamma^t(x, \cdot) - \pi(\cdot) } = 0,
\ee
and satisfies the
minorization condition
\be
\label{eq:minor}
K_\gamma(x,\cdot) \ge m_0 \xi(\cdot),
\quad \forall x \in \R,
\ee
where $\xi(\cdot)$ is a probability measure on $(\R,\sigalg)$ 
%that has a
%density with respect to $\lambda,$
and $m_0$ is some positive constant.
\end{assumption}

%\begin{assumption}
%\label{assump:3}
%The kernels $K_\gamma(\cdot,\cdot)$ are total-variation
%continuous with respect to $\gamma$, uniformly in their first arguments.
%That is, for every $\epsilon>0$, there exists $\delta(\epsilon)$ such
%that
%\be
%|\gamma_1 - \gamma_2| < \delta(\epsilon)
%\rightarrow
%\TV{ K_{\gamma_1}(x, \cdot) - K_{\gamma_2}(x, \cdot) }
% < \epsilon, \quad \forall x \in \R.
%\ee
%\end{assumption}

Assumptions~\ref{assump:1} and \ref{assump:2}
%, and \ref{assumpt:3}
are not restrictive.  The first can be satisfied by requiring that
$\abs{ \Gamma_{t+1} - \Gamma_t } = o(t^{-\alpha})$ for some $\alpha>1$.
Intuitively, this is just a form of what \cite{RobertsRosenthal} refer to 
as ``diminishing adaptation''.
Assumption~\ref{assump:2} ensures that all possible adaptive kernels
mix at a minimal rate, and that any composition of kernels also mixes
at that rate. The latter follows from the well-known bound on the
contraction coefficient $\theta$ via the minorization constant $m_0$
in (\ref{eq:minor}):
\beq
\theta \leq 1-m_0;
\eeq
see, for example, Lemma 2.2.3 in 
\citet{kont07-thesis}.
One way to construct such a family of kernels is to choose a family of Metropolis-Hastings kernels in which the proposal distributions all share a common component which does not depend on the current state.

Assumption~\ref{assump:4} 
is a natural Feller-type condition; in particular, 
it is satisfied by most Metropolis-Hastings chains, and is fairly easily checked in practice.
%it is certainly satisfied by a 
%Metropolis-Hastings chain.
\hide{
%is 
seems
more restrictive.  
It requires the chain to have a density with respect to a product measure.  This is generally easy to establish for discrete space Metropolis-Hastings chains (e.g. those defined on the integers), as these chains typically have densities with respect to a product counting measure.  
On the other hand, a Metropolis-Hastings chain taking values on the real line generally has a measure which assigns positive probability to events of the form $\{X_0=X_1=X_2\}$, for which there is no
%and it is not generally clear how one would find a 
dominating product measure.
However, our methods easily extend to this case. Recall the rapid convergence result (\ref{eq:expfast}) in Theorem \ref{thm:lln}.
}

Using our main result along with these assumptions, we are in a position to state conditions under which 
Theorem 23 of \cite{RobertsRosenthal} can be strengthened to establish strong instead of weak convergence.

\bethn
\label{thm:adaptslln}
Suppose that an adaptive Markov chain satisfies Assumptions~\ref{assump:1},
\ref{assump:2}, and \ref{assump:4}.  Then we have
$$
\hat{P}_n(A) \rightarrow \pi(A), 
\text{ a.s.},
$$
for each $A \in {\cal F}$, where $\hat{P}_n(\cdot)$ is the empirical measure defined in (\ref{eq:pndef}).
\enthn

\begin{proof}
The minorization condition in Assumption~\ref{assump:2} ensures that the simultaneous uniform ergodicity condition of Theorem 5 of \cite{RobertsRosenthal} holds.  
Assumption~\ref{assump:1} ensures that the diminishing adaption condition of the same theorem is also satisfied.  
Then we have
\be
%\E[ \abs{\hat{P}_n(A) - \pi(A)} ] 
\abs{\E \hat{P}_n(A) - \pi(A) }
\le 
  \epsilon + \pr{\abs{\hat{P}_n(A)- \pi(A)} > \epsilon },
\ee
and the latter probability is bounded by a quantity independent of $A$, 
%so
%\be
%\TV{ \E\hat{P}_n(\cdot) - \pi(\cdot) }
%\le
%\epsilon + \sup_A \pr{\abs{\hat{P}_n(A)- \pi(A)} > \epsilon},
%\ee
%but the supremum in the last expression is bounded uniformly by a finite constant 
which goes to zero, 
by Theorem 23 of \cite{RobertsRosenthal} (the uniformity is not stated explicitly in the Theorem but is established in the proof that they give).  
This establishes the uniformly converging expectation condition of our Theorem~\ref{thm:lln}.   
The minorization condition 
(\ref{eq:minor})
ensures that the mixing coefficients $\bar\eta_{ij}$ decay as $(1-m_0)^{j-i}$.
%{\bf Also mixing matrix is well behaved because of minorization condition}.  
%Finally, the existence of a density with respect to a product measure comes from Assumption~\ref{assump:4}.
To establish almost-sure convergence, we are going to argue along the lines of 
Theorem~\ref{thm:lln} (the latter is not applicable directly, since
%We are almost ready to apply
%%Thus the conditions of 
%Theorem~\ref{thm:lln},
%the only possible glitch being that 
the adaptive Markov chain measure $P$ might not have a density with respect to the Lebesgue measure).
Let $n_0=n_0(\eps)$ 
be such that
$\abs{\E \hat{P}_n(A) - \pi(A) } < \eps$
for all $n>n_0(\eps)$. As in Theorem \ref{thm:lln}, the quantity we want to bound is
\beqn
\label{eq:contprob}
P\paren{ \abs{ \hat{P}_n(A) - \pi(A) } > t+\eps}
\eeqn
for $\eps,t>0$ and $n>n_0(\eps)$.

We now approximate the chain $P$ by a finite-state chain $P'$ 
%We get around this by 
%%fixing a set $A\in\calF$ and 
%approximating the chain $P$ 
%the a finite-state chain $\tilde P$ 
induced
by finite partitions, as shown in the Appendix. 
%Since 
%Now
%Using arguments similar to those 
It follows from the arguments
in the Appendix
that
for any $E\in\calF^n$
we can find a partition of the 
state space
and a finite-state chain $P'$ on $(\Omega',\calF')^n$ so that
$P'(\tilde E)$
is
arbitrarily close to $P(E)$
for some 
$\tilde E\in(\calF')^n$.
Therefore, 
%it follows that 
$\hat P_n(A)$ and $\pi(A)$ can be made
arbitrarily close to their finite-state analogues
$\hat P'_n(\tilde A)$ and $\pi'(\tilde A)$
and in particular,
\beqn
\label{eq:discrprob}
P'\paren{ \abs{ \hat{P}'_n(\tilde A) - \pi'(\tilde A) } > t+\eps}
\eeqn
approximates the expression in (\ref{eq:contprob}).
%[Note that we are abusing notation by writing $P'(E)$ as opposed to 
%What's more, 
Furthermore, 
Lemma \ref{lem:cont2fin} shows that for a sufficiently refined partition, the finite-state chain $P'$
will have mixing coefficients
arbitrarily close
to those of $P$.
%approximate the  of $P$ 
%arbitrarily closely.
To conclude the proof it suffices to apply Theorem~\ref{thm:lln} and Corollary \ref{cor:asconv}.
%For a fixed $A\in\calF$, we may approximate the probabilities and mixing coefficients 
%Theorem~\ref{thm:lln} applies to any finite-state chain.
%%are met, and the 
%The
%strong law follows directly from an application of the Borel-Cantelli lemma.
\end{proof}

\section*{Acknowledgements}
We would like to thank Ofer Zeitouni for helpful discussions.

%\section*{Appendix}
\begin{appendix}

\section{Extending finite-state inequalities to more general spaces}

\newcommand{\hX}{\hiddenX}
\newcommand{\oX}{\observedX}
\newcommand{\ohX}{\calW}
\newcommand{\ohx}{x^{\osymb\hsymb}}

In bounding the mixing coefficients, measure-theoretic technicalities
tend to play a peripheral role. Indeed, the $\X=\set{0,1}$ case
already captures most of the proof complexity. The mixing results we
proved for finite $\X$ extend verbatim to $\X=\N$, and under mild
%topological
continuity
assumptions, to much more general measure spaces. 

The 
inequality
%concentration result 
in Theorem \ref{thm:mgalebd} applies to a broad class of state spaces, including $\X=\R$. The latter
result estabishes that the $\eta$-mixing coefficients of a random process control the 
%probability of 
concentration of Lipschitz path functionals about their means. This immediately implied the strong law of large numbers
in Section \ref{sec:slln}.

In Theorem \ref{thm:mmconc} we showed how to control the $\eta_{ij}$ in terms if the Markov contraction coefficients; this was done for finite
state spaces. In this Appendix, we extend these bounds to the continuous case.

\subsection{Markov marginal chains}

In the continuous case we shall consider the ``observed'' state space $\oX$
and the ``hidden'' state space $\hX$, with the total space
$\ohX=\oX\times\hX$; one may take $\oX=\hX=\R$. 
We take the usual Borel $\sigma$-algebra on 
$\ohX$, which is a product of the corresponding
$\sigma$-algebras on $\oX$ and $\hX$.
%We will assume a positive Borel measure $\mu$ on $\ohX$,
%assuming further that the $\sigma$-algebra on $\ohX$ is a product of
%$\sigma$-algebras on $\oX$ and $\hX$, and $\mu=\mu_{\oX}\tp\mu_{\hX}$.
%\hide{
%We will assume that
%$\ohX$ is endowed with a metric $\rho$, which induces a Borel
%$\sigma$-algebra, and the latter is equipped with a positive measure $\mu$.
%}
%%, $\hiddenX=\observedX=\R$,
%%[perhaps write $\calX,\calY$ instead of $\hiddenX,\observedX$ ??]
%%and $\mu$ will denote the Lebesgue measure on $\R$.
An MMC
%$\sseq{X}{1}{n}$
over $\oX^n$
is
obtained by first defining the Markov process $W_t=(X_t,Y_t)$, with
values in
$\ohX
%=\observedX\times\hiddenX
$,
induced by the kernels
%defined by the Markov kernels
$\set{K_i}_{0\leq i<n}$.
%from $\ohX=\observedX\times\hiddenX$ to itself.
Thus if $A_i\subset\ohX$, $1\leq i\leq n$ are measurable, then
\beq
\pr{W_1\in A_1,\ldots,W_n\in A_n}
&=&
%\int_{\ohX^n}d\ohx_1\ldots d\ohx_n
\int_{A_1}\ldots \int_{A_n}
\prod_{i=1}^n
%K_i(\ohx_{i-1},d\ohx_i)
K_i(w_{i-1},d w_i).
\eeq
We will write $K_0$ for the initial distribution of the MMC, and
expressions such as $K_0(x_0,d x_1)$ are to be interpreted as
$K_0(dx_1)$. 
%We also assume henceforth that the Markov measure $P$ on $\ohX^n$ satisfies
%$P\abscont\mu^n$.
Note that we are not assuming anything about the density of $P$; the latter may well not exist with respect to the 
Lebesgue (or any other product) measure on $\ohX^n$.

As in the discrete case, the MMC is obtained by marginalizing out the
``hidden'' component $Y$:
\beq
\pr{X_1\in B_1,\ldots,X_n\in B_n}
&=&
\int_{\hiddenX^n}
\int_{B_1}\ldots \int_{B_n}
\prod_{i=1}^n
%K_i((\o x_{i-1},\h x_{i-1}),(d\o x_{i},d\h x_{i}))
K_i((x_{i-1},y_{i-1}),(d x_{i},d y_{i}))
\eeq
for measurable sets $B_i\subset\observedX$.

The definition of the contraction coefficient readily generalizes to the continuous case:
\beqn
\label{eq:thadefR}
\tha_i &=&
%\essup
\sup
_{w,w'\in \ohX}
\TV{
K_i\paren{w,\cdot} -
K_i\paren{w',\cdot}
}.
\eeqn
We will use the notation $\tha_i(P)$ if we wish to make explicit the dependence on the particular Markov measure.

Similarly,
we define
\beq
\bar\eta_{ij} &=&
%\essup
\sup
_{y\in\ohX^{i-1},w,w'\in\ohX}
\eta_{ij}(y,w,w'),
\eeq
and write $\bar\eta_{ij}(P)$ to make explicit the particular measure.

All sets below are assumed to be Borel-measurable.
Recall that a 
%(finite) 
{\em partition} of a set $E$ is a collection of disjoint sets whose union is $E$.
%Any 
Whenever
$\calA$
is a collection of subsets of $E$,
$\calA$ induces an equivalence relation on $E$ as follows:
\beq
x\equiv_\calA y &\Longleftrightarrow&
\pred{x\in A}=\pred{y\in A} \text{ for all } x,y\in E, A\in\calA 
%x,y\in A \text{ for some } A\in\calA
.
\eeq
Such an equivalence relation in turn induces a partition on $E$ (whose members are the equivalence classes);
we shall call such a partition the {\em refining partition} of $\calA$. Notice that if $\calA$ is finite
then so is its refining partition.
%we say that a partition $\calU$
%{\em refines} $\calA$ if for all $A,B\in\calA$, $A\cap B$ is the union of members of $\calU$.
%For example, the algebra of sets generated by a finite collection $\calA$
%contains a 

Let $P$ be a Markov chain on $\ohX^n$, as above. Any finite 
%(and measurable) 
partition 
$\set{V_k:1\leq k\leq m}$ of $\ohX$ induces a Markov chain $\hat P$ 
with kernels $\set{\hat K_i}_{0\leq i<n}$
%on 
%$\hat\ohX$, 
%where
over
the finite state space
$\hat\ohX=\set{1,\ldots,m}$, as follows:
\beq
\hat K_i(k', k) &=& P(W_{i+1} \in V_k \gn W_i \in V_{k'}).
\eeq
If $\ohX=\oX\times\hX$, $\set{S_k}$ is a partition of $\oX$, and 
$\set{T_\ell}$ is a partition of $\hX$, then these two partitions induce a partition on
$\ohX$ in the obvious way. 
We will write $\hat \ohX=\hat \oX\times\hat\hX$ to denote the state spaces obtained by identifying
the partition blocks with states.
%Furthermore, identifying the partition elements as states, 

\belen
\label{lem:cont2fin}
Let $P$ be a Markov chain on 
$\ohX^n=(\oX\times\hX)^n$
%, with $\tha_i(P)=t$. 
and $Q$ be the induced MMC on $\oX^n$.
Then, 
assuming that the kernel generating $P$ satisfies the Feller
%-type 
continuity condition
in Assumption~\ref{assump:4}, we have
that
for any
$\eps>0$ there are finite partitions of 
%$\ohX$ 
$\oX$ and $\hX$
such that the induced Markov chain $\hat P$ 
on $\hat\ohX^n=(\hat\oX\times\hat\hX)^n$
and MMC $\hat Q$ 
on $\hat\oX^n$
satisfy
\beqn
\label{eq:thapart}
%\tha_i(P)-\eps \;\leq\; \tha_i(\hat P) \; \leq\; \tha_i(P)
\abs{\tha_i(P)-\tha_i(\hat P)} &<& \eps
\eeqn
and
\beqn
\label{eq:etapart}
\abs{\bar\eta_{ij}(Q)-\bar\eta_{ij}(\hat Q)} &<& \eps
%\bar\eta_{ij}(Q)-\eps \;\leq\; \bar\eta_{ij}(\hat Q) \; \leq\; \bar\eta_{ij}(Q)
,
\eeqn
for $1\leq i<j\leq n$.
\enlen
\bepf
Fix an $\eps>0$.
Let us construct the requisite partition for (\ref{eq:thapart}). 
Fix an $1\leq i<n$ and
%suppose
let
$t=\tha_i(P)$. Then 
by definition of $\TV{\cdot}$
there are 
$w%=(x,y)
$ 
and 
$w'%=(x',y')
$ in $\ohX$,
and an $E\subset\ohX$
%and a measurable, disjoint
%$A_k\subset \oX$, $k=1,\ldots,m_A$,
%$B_\ell\subset \hX$, $\ell=1,\ldots,m_B$,
%% $A\subset\ohX$
such that
\beq
t-\eps/2 \;\leq\; K_i(w,E)- K_i(w',E) \; \leq\; t
%+\eps/2
.
\eeq
Furthermore, $E$ may be approximated by a finite union of rectangles
$A_k \times B_\ell$, with 
$A_k\subset \oX$, $k=1,\ldots,m_A$,
and
$B_\ell\subset \hX$, $\ell=1,\ldots,m_B$:
%where $E\subset\ohX$ is given by
\beqn
\label{eq:Edecomp}
\tilde E &=& \bigcup_{k=1,\ldots,m_A} \bigcup_{\ell=1,\ldots,m_B} A_k \times B_\ell,
\eeqn
so that
\beq
t-\eps \;\leq\; K_i(w,\tilde E)- K_i(w',\tilde E) \; \leq\; t
%+\eps/2
.
\eeq
By Feller continuity, we have 
\beq
K_i(w,E) &=& \lim_{
\alpha\to\infty
}
P(W_{i+1}\in E \gn W_i \in U_\alpha)
\eeq
where $U_1\supset U_2\supset \ldots \ni w$ is a sequence of neighborhoods
{\em shrinking to $w$} in the sense that $\cap_\alpha U_\alpha=w$.
(We are assuming without loss of generality that $P(U_\alpha)>0$ for
all $\alpha$.)
%; if no such sequence can be found, the conditional
%distribution $K_i(w,\cdot)$ may be (re)defined arbitrarily.]

Therefore, taking $w=(x,y)$ and $w'=(x',y')$, we have that there are neighborhoods 
$C,C' \subset\oX$
of 
$x$ and $x'$ (respectively),
as well as analogous neighborhoods
$D,D'\subset\hX$ 
of $y$ and $y'$, such that
$$
P(W_{i+1}\in E \gn W_i \in C\times D)
-
P(W_{i+1}\in E \gn W_i \in C'\times D')
$$
is arbitrarily close to $K_i(w,\tilde E)- K_i(w',\tilde E)$.

%Let let $C,C'\subset\oX$ be sets of 
%small but
%positive $\mu_{\oX}$ measure so that $x\in C$ and $x'\in C'$
%and define $D,D'\subset\hX$ analogously to contain $y$ and $y'$, respectively.
Let $\set{S_a}$ be the partition of $\oX$ refining the sets $\set{A_k}$, $C$, and $C'$.
Similarly, let
$\set{T_b}$ be the partition of $\hX$ refining the sets $\set{B_\ell}$, $D$, and $D'$.
Identify
the partition blocks with states
and induce the finite-state Markov chain $\hat P$ on
$\hat \ohX=\hat \oX\times\hat\hX$.
By the approximation argument above,
%Since
%\beq
%K_i(w,E) &=& \lim_{
%\mu_{\oX}(C)\to0,\mu_{\hX}(D)\to0
%}
%P(W_{i+1}\in E \gn W_i \in C\times D)
%,
%\eeq
%as $\mu_{\oX}(C)\to0$ and $\mu_{\hX}(D)\to0$, 
%it is straightforward to verify that
there exist $C,C',D,D'$ such that
$\hat P$ satisfies (\ref{eq:thapart}) for the given $i$. Repeating this process for each $i=1,
\ldots,n-1$ and picking a finite partition of $\oX$ (respectively, $\hX$) 
that simultaneously refines the partitions for each $i$, 
we establish (\ref{eq:thapart}) for all $i$.
%Let us denote the $\oX$-partition thus obtained by $\set{S^*_a}$
%and the corresponding $\hX$-partition by $\set{T^*_b}$.

Now we turn to (\ref{eq:etapart}). The basic technique is the same as the one used to
show (\ref{eq:thapart}).
%, but the notation is bound to get messier.
Fix $1\leq i<j\leq n$ and
let
$h=\bar\eta_{ij}(Q)$.
Then there are 
$\sseq{x}{1}{i-1}\in\oX^{i-1}$, $x_i,x_i'\in\oX$
and 
an $A\subset\oX^{n-j+1}$ such that
%a measurable $A=A_{j}\times A_{j+1}\times\ldots\times A_n \subset\oX^{n-j+1}$
%such that\footnote{
%This is not entirely accurate; technically, we should decompose each 
%$A_k$
%as a finite partition 
%%before
%$A_k=\cup_{\ell}A\supr{k}_\ell$,
%move the $\cup$ to the outside and the $\times$ to the inside,
%%taking the product, 
%as done with $E$ 
%%above
%in (\ref{eq:Edecomp}).}
\beq
h-\eps/2 \;\leq\; 
Q(\sseq{X}{j}{n}\in  A \gn \sseq{X}{1}{i} = \sseq{x}{1}{i} )
-
Q(\sseq{X}{j}{n}\in A \gn \sseq{X}{1}{i} = \sseq{x}{1}{i-1} x')
\; \leq\; h.
\eeq

As done with $E$ and $\tilde E$ above, 
$Q(\sseq{X}{j}{n}\in  A \gn \sseq{X}{1}{i} = \sseq{x}{1}{i} )$
may be approximated by
$Q(\sseq{X}{j}{n}\in  \tilde A \gn \sseq{X}{1}{i} = \sseq{x}{1}{i} )$
where $\tilde A$ is now an $(n-j+1)$-fold product of expressions of the type 
$\cup_{k} A_k$, with $A_k\subset \oX$.

Again as above, we take small neighborhoods around $x_1,x_2,\ldots,x_{i},x'_{i}$ to
%appearing in the r.h.s. of
%(\ref{eq:Edecomp})
%
%Letting $x_i'\in A_0\subset\oX$, $x_k\in A_k\subset\oX$, $k=1,\ldots,i$ be 
%small neighborhoods
%%sets of 
%%small $\mu_{\oX}$-measure 
%and taking
%%$\set{S^{**}_a}$ to be 
%the partition of $\oX$ that 
%%refines $\set{S^{*}_a}$ and 
%the $\set{A_k:0\leq k\leq n}$, we 
obtain a partition which satisfies 
(\ref{eq:etapart}) for the given $i,j$. Refining such partitions 
simultaneously for
$1\leq i<j\leq n$, we establish (\ref{eq:etapart}) for all these values.

Finally, the partitions obtained in the course of proving
(\ref{eq:thapart}) and (\ref{eq:etapart}) can once again be refined to make the two inequalities
hold simultaneously.
\enpf
\becon
\label{cor:mainresmmc}
Let $P$ be a Markov chain on 
$\ohX^n=(\oX\times\hX)^n$
and $Q$ be the induced MMC on $\oX^n$. Then
\beq
\bar\eta_{ij}(Q)
&\leq&
\prod_{k=i}^{j-1} \tha_k(P).
\eeq
\encon
\bepf
Immediate consequence of the corresponding claim for finite $\oX$ and $\hX$.
\enpf

\subsection{Tensorization}

\belen
\label{lem:tensR}
Suppose the Markov chain $P$ on $(\oX\times\hX)^n$ 
is defined by transition
kernels $\set{K_i}_{0\leq i\leq n}$, which have the
following special structure:
\beqn
\label{eq:kdecompR}
K_{i}\paren{(x',y'),(x,y)}
&=&
A_i(x \gn (x',y'))
B_i(y \gn (x',y')).
\eeqn
Then we have
\beqn
\label{eq:thabR}
\tha_i &\leq& \alpha_i+\beta_i-\alpha_i\beta_i
\eeqn
where 
%$\tha_i$ is defined in (\ref{eq:thadef}), and
\beq
\alpha_i &=& 
\sup_{(x',y'),(x'',y'')\in \oX\times\hX}
\TV{
A_i(\cdot \gn (x',y')) -
A_i(\cdot \gn (x'',y'')) 
},\\
\beta_i &=& 
\sup_{(x',y'),(x'',y'')\in \oX\times\hX}
\TV{
B_i(\cdot \gn (x',y')) -
B_i(\cdot \gn (x'',y''))
}.
\eeq
\enlen
\bepf
The same technique of approximating $P$ by a finite-state Markov chain $\hat P$ as
employed in Lemma \ref{lem:cont2fin} may be used here; details are omitted.
\enpf

\subsection{Adaptive Markov chains}
Let $\Gamma$ be an index set 
%endowed with a Borel measure $\mu_\Gamma$
and $\set{K_\g}_{\g\in\Gamma}$ be a collection of
Markov transition kernels, $K_\g:\X\to\X$. 
%We assume a measure $\mu_\X$ on $\X$ and a product measure $\mu=\mu_\X\tp\mu_\Gamma$.
For a given sequence
$\sseq{\g}{1}{n-1}\in\Gamma^{n-1}$, we have a Markov measure on
$\X^n$:
\beq
dP_{\sseq{\g}{1}{n-1}}(x) &=& 
\prod_{i=0}^{n-1}K_{\g_i}(x_{i}, dx_{i+1}),
\qquad
x\in\X^n.
\eeq
For $1\leq i<n$, $x\in\X$, and $\g'\in\Gamma$, let $g_i(\cdot\gn x,\g')$ be a probability
measure on $\Gamma$. 
Together,
$\set{K_\gamma}$ and $\set{g_i}$ define a measure $Q$ on $\X^n$, which we'll call
an {\em adaptive Markov} measure:
\beqn
\label{eq:adapt-defR}
dQ(x) &=& 
\int_{\Gamma^{n}}
%d\mu_\Gamma(\sseq{\g}{1}{n})
%%\sum_{\sseq{\g}{1}{n}\in\Gamma^{n}}
%%p_0(x_1,\g_1)
\prod_{i=1}^{n-1}[g_i(d\g_{i+1}\gn x_i,\g_{i})K_{\g_i}(x_{i}, dx_{i+1})],
\eeqn
where $\g_1$ is a dummy index to make $g_1(\cdot\gn\cdot)$ well-defined.

Define the following contraction coefficients

\beqn
\label{eq:kapdefR}
\kappa
&=& \sup_{w,w'\in\X,\g,\g'\in\Gamma}\TV{K_{\g}(w,\cdot)-K_{\g'}(w',\cdot)}.
\eeqn

and

\beqn
\label{eq:lamdefR}
\lambda_i
&=& 
%\sup_{1\leq i<n} 
\sup_{w,w'\in\X,\g,\g'\in\Gamma}
\TV{g_i(\cdot\gn w,\g)-g_i(\cdot\gn w',\g')}.
\eeqn

Then
\bethn
For an adaptive Markov measure $Q$ on $\X^n$
we have
\beq
\bar\eta_{ij} &\leq& \tha_i \tha_{i+1}\cdots \tha_{j-1}
\eeq
where 
$$ \tha_i = \kappa+\lambda_i - \kappa\lambda_i.$$
\enthn
\bepf
First we observe that the
adaptive Markov process 
defined in (\ref{eq:adapt-defR})
is in fact a MMC,
with $\oX=\X$ and $\hX=\Gamma$;
thus 
%Theorem \ref{thm:pomc} 
Corollary \ref{cor:mainresmmc}
applies.
Furthermore, 
%its 
the MMC
kernel 
(denoted here by $P_i$ instead of $K_i$ to avoid confusion with $K_\g$)
decomposes as in (\ref{eq:kdecompR}):
\beqn
\label{eq:pgR}
P_{i}\paren{
(x',\g'),(dx,d\g)}
&=&
g_i(d\g \gn x',\g')
K_{\g'}(x', dx).
\eeqn
The claim is proved by applying Lemma \ref{lem:tensR} to (\ref{eq:pgR}).
\enpf

\end{appendix}

\bibliography{slln}

\end{document}